\documentclass[A4paper,11pt]{article}
\usepackage{amsfonts}
\usepackage{amsmath}
\usepackage{amssymb}
\usepackage{graphicx}
\usepackage{amscd}
\begin{document}
\date{}
\title{Homology of torus links
\author{Marko Sto\v si\'c \\
Instituto de Sistemas e Rob\'otica and CAMGSD,\\
Instituto Superior T\'ecnico, 
Av. Rovisco Pais 1\\
1049-001 Lisbon, 
Portugal\\
e-mail: mstosic@math.ist.utl.pt
}
}

\newtheorem{theorem}{Theorem}
\newtheorem{acknowledgment}[theorem]{Acknowledgment}
\newtheorem{algorithm}[theorem]{Algorithm}
\newtheorem{axiom}[theorem]{axiom}
\newtheorem{case}[theorem]{Case}
\newtheorem{claim}[theorem]{Claim}
\newtheorem{conclusion}[theorem]{Conclusion}
\newtheorem{condition}[theorem]{Condition}
\newtheorem{conjecture}[theorem]{Conjecture}
\newtheorem{corollary}[theorem]{Corollary}
\newtheorem{criterion}[theorem]{Criterion}
\newtheorem{definition}{Def\mbox{}inition}
\newtheorem{example}{Example}
\newtheorem{exercise}[theorem]{Exercise}
\newtheorem{lemma}{\indent Lemma}
\newtheorem{notation}[theorem]{Notation}
\newtheorem{problem}[theorem]{Problem}
\newtheorem{proposition}{Proposition}
\newtheorem{remark}[theorem]{Remark}
\newtheorem{solution}[theorem]{Solution}
\newtheorem{summary}[theorem]{Summary}
\newcommand{\ud}{\mathrm{d}}

\def\gcd{\mathop{\rm gcd}}
\def\Ker{\mathop{\rm Ker}}
\def\max{\mathop{\rm max}}
\def\map{\mathop{\rm map}}
\def\lcm{\mathop{\rm lcm}}
\def\kraj{\hfill\rule{6pt}{6pt}}
\def\diag{\mathop{\rm diag}}
\def\span{\mathop{\rm span}}
\def\deg{\mathop{\rm deg}}
\def\rank{\mathop{\rm rank}}
\def\sgn{\mathop{\rm sgn}}
\def\kvn{\{n\}_q}
\def\F{\mathbb{F}}
\def\R{\mathbb{R}}
\def\C{\mathcal{C}}
\def\P{\mathcal{P}}
\def\A{\mathcal{A}}
\def\D{\mathcal{D}}
\def\H{\mathcal{H}}
\def\N{\mathbb{N}}
\def\K{\mathbb{K}}
\def\Z{\mathbb{Z}}
\def\Q{\mathbb{Q}}
\def\S{\mathbb{S}}
\def\X{\qbezier(0.00,0.00)(0.50,1.00)(1.00,2.00)
\qbezier(1.00,0.00)(0.80,0.40)(0.60,0.80)
\qbezier(0.00,2.00)(0.20,1.60)(0.40,1.20)
}

\def\Y{\qbezier(1.00,0.00)(0.50,1.00)(0.00,2.00)
\qbezier(0.00,0.00)(0.20,0.40)(0.40,0.80)
\qbezier(1.00,2.00)(0.80,1.60)(0.60,1.20)
}

\def\O{\qbezier(0.00,0.00)(0.20,1.00)(0.00,2.00)
\qbezier(1.00,0.00)(0.80,1.00)(1.00,2.00)
}
\def\OF{\qbezier(0.00,0.00)(0.00,1.00)(0.00,2.00)
\qbezier(1.00,0.00)(1.00,1.00)(1.00,2.00)
}

\def\l{
\qbezier(0.00,0.00)(0.50,1.20)(1.00,0.00)
\qbezier(0.00,2.00)(0.50,0.80)(1.00,2.00)
}

\def\LPP{
\qbezier(0.00,0.00)(0.35,0.70)(0.70,1.40)
\qbezier(0.70,0.60)(0.65,0.70)(0.60,0.80)
\qbezier(0.00,2.00)(0.20,1.60)(0.40,1.20)
\qbezier(0.7,1.4)(1.3,2.6)(1.3,1)
\qbezier(0.7,0.6)(1.3,-0.6)(1.3,1)}

\def\LPM{\qbezier(0.00,0.00)(0.2,0.40)(0.40,0.80)
\qbezier(0.60,1.20)(0.65,1.30)(0.70,1.40)
\qbezier(0.00,2.00)(0.35,1.30)(0.70,0.60)
\qbezier(0.7,1.4)(1.3,2.6)(1.3,1)
\qbezier(0.7,0.6)(1.3,-0.6)(1.3,1)}

\def\XS{\qbezier(0.00,0.00)(0.50,1.00)(1.00,2.00)
\qbezier(1.00,0.00)(0.80,0.40)(0.60,0.80)
\qbezier(0.00,2.00)(0.20,1.60)(0.40,1.20)

\qbezier(0.0,2.00)(0.0,1.85)(0,1.70)
\qbezier(1.00,2.00)(1,1.85)(1.0,1.70)
\qbezier(0.00,2.00)(0.10,1.90)(0.20,1.80)
\qbezier(1.00,2.00)(0.9,1.9)(0.8,1.80)}

\def\YS{\qbezier(1.00,0.00)(0.50,1.00)(0.00,2.00)
\qbezier(0.00,0.00)(0.20,0.40)(0.40,0.80)
\qbezier(1.00,2.00)(0.80,1.60)(0.60,1.2)
\qbezier(0.0,2.00)(0.0,1.85)(0,1.70)
\qbezier(1.00,2.00)(1,1.85)(1.0,1.70)
\qbezier(0.00,2.00)(0.10,1.90)(0.20,1.80)
\qbezier(1.00,2.00)(0.9,1.9)(0.8,1.80)}

\def\GR{\qbezier(0.00,0.00)(0.25,0.25)(0.50,0.50)
\qbezier(1.00,0.00)(0.75,0.25)(0.50,0.50)
\qbezier(0.00,2.00)(0.25,1.75)(0.50,1.5)
\qbezier(1.0,2.00)(0.75,1.75)(0.5,1.50)

\qbezier(0.00,2.00)(0,1.9)(0,1.80)
\qbezier(0.00,2.00)(0.1,2)(0.20,2)
\qbezier(1.00,2.00)(1,1.9)(1,1.80)
\qbezier(1.00,2.00)(0.9,2)(0.8,2)
\linethickness{2.5pt}
\qbezier(0.5,0.5)(0.5,1)(0.50,1.5)}

\def\GRF{\qbezier(0.00,0.00)(0.25,0.25)(0.50,0.50)
\qbezier(1.00,0.00)(0.75,0.25)(0.50,0.50)
\qbezier(0.00,2.00)(0.25,1.75)(0.50,1.5)
\qbezier(1.0,2.00)(0.75,1.75)(0.5,1.50)

\linethickness{2.5pt}
\qbezier(0.5,0.5)(0.5,1)(0.50,1.5)}

\def\OS{\qbezier(0.00,0.00)(0.60,1.00)(0.00,2.00)
\qbezier(1.00,0.00)(0.40,1.00)(1.00,2.00)
\qbezier(0.0,2.00)(0.0,1.85)(0,1.70)
\qbezier(1.00,2.00)(1,1.85)(1.0,1.70)
\qbezier(0.00,2.00)(0.10,1.90)(0.20,1.80)
\qbezier(1.00,2.00)(0.9,1.9)(0.8,1.80)}

\maketitle

\begin{abstract}
In this paper we show that there is a cut-off in the Khovanov homology of $(2k,2kn)$-torus links, namely that the maximal homological degree of non-zero homology group of $(2k,2kn)$-torus link is  $2k^2n$. Furthermore, we calculate explicitely the homology groups in homological degree $2k^2n$ and prove that it coincides with the centre of the ring $H^k$ of crossingless matchings, introduced by M. Khovanov in \cite{tan}. Also we give an explicit formula for the ranks of the homology groups of $(3,n)$-torus knots for every $n\in\N$.
\end{abstract}

\section{Introduction}

The Khovanov homology of links was introduced by M. Khovanov  in \cite{kov}. For every link, the graded chain complex is defined whose homology groups are link invariants, and whose graded Euler characteristic is equal to the Jones polynomial of a link.\\
\indent One of the main advantages of the theory is that its definition is combinatorial and there is a straightforward algorithm for computing the homology groups. Consequently, there are various computer programs \cite{bnprog}, \cite{shum} that  efficiently calculate  homology groups for links with 50 crossings and more, and the results obtained are the basis for many conjectures about the form of homology. 
However, not much is proved about the homology of  large classes of knots. It is known (see \cite{ap},\cite{lee},\cite{lee2}) that ranks of the homology groups of alternating knots are determined by the Jones polynomial and knot signature, and that the homology of alternating knots is  contained in two diagonals.\\
\indent On the other side is the class of torus knots, which is ``the most non-alternating" class of knots, and  not much is known about their homology.  By using the computer programs, it was found that the homology of torus knots occupies lots of diagonals, and that it is rather complicated. Also, the torus knots are the smallest knots (in the sense of the number of crossings) where the properties, like homological thickness and torsions of high order appear (\cite{bnprog}, \cite{shtor}).
 
In our previous papers (\cite{mojtorus}, \cite{teza}), we have already obtained various properties of the homology of torus knots $T_{p,q}$. More precisely, we proved that the homology groups $H^i(D_{p,q})$ and $H^i(D_{p,q-1})$ coincide up to certain homological degree $i$ (at least $p+q-3$), and hence that there exists stable homology of torus knots when $q$ tends to infinity. This also  implies that the non-alternating torus knots are homologically thick, and that the only almost alternating torus knots are $T_{3,4}$ and $T_{3,5}$, thus solving a conjecture from \cite{adams}.

The basic ingredient that we are using is the long exact sequence in homology, which is particularly powerful in the case of torus knots. Here we use similar approach to calculate ``the upper" part (large homological degrees) of the homology of torus links.\\

\indent In this paper we prove further properties for torus links with even number of components.  Namely we show that torus link $T_{2k,2kn}$ has nonzero homology only up to the homological degree $2k^2n$. Since we assumed that torus knots $T_{p,q}$ are positive and the number of crossings of $D_{2k,2kn}$ is $(2k-1)2kn$, this means that we have proved that the homology of a torus link is contained only in the left half of the plane. 

\indent This can be seen even better if we observe torus link $T'_{2k,2kn}$  which coincides with $T_{2k,2kn}$ when exactly $k$ components are with reversed orientation. Then this torus link has $2k^2n$ negative crossings and  $2k(k-1)n$ positive crossings, while its homology coincides with the appropriately shifted homology of $T_{2k,2kn}$. In other words, the statement from the previous paragraph  says that the homology of $T'_{2k,2kn}$ is trivial in the positive homological degrees. \\
\indent Analogously as the previous result, we also obtained that the similar fact is true when we observe the $q$-gradings. Namely we show that the $q$-gradings of non-zero homology groups of $T_{2k,2kn}$ are less or equal than $6k^2n$. In terms of torus link $T'_{2k,2kn}$, this means that the $q$-gradings of non-zero homology groups are non-positive, and hence the homology of torus link $T'_{2k,2kn}$ is contained in the lower-left quadrant.

Analogously as in the  case of torus links, for a general torus knot $T_{p,q}$ we obtain that it has trivial homology for homological degrees bigger than $pq/2$.\\

\indent In \cite{prz}, J. Przytycki has shown that the limit of Khovanov homology of torus links $T'_{2,2n}$ when $n$ tends to infinity, coincides with the Hochschild cohomology of the ring $H^1$ (isomorphic to $\Z[X]/(X^2)$, see \cite{tan}, \cite{kovspr}). Recently, it was conjectured by M. Khovanov and L. Rozansky in \cite{KRS2} that there exists limit of the homology of torus links $T'_{2k,2kn}$ when $n$ tends to infinity and that it coincides with the Hochschild cohomology of the ring $H^k$ of crossingless matchings. In particular, this means that the zeroth homology group of $T'_{2k,2kn}$ coincides with the zeroth Hochschild cohomology group of $H^k$, which equals the center of $H^k$.

In this paper we prove the latter fact: first we calculate explicitely the homology groups of $T'_{2k,2kn}$ in the homological degree 0, then we determine the center of the ring $H^k$, and finally we conclude that they coincide, up to appropriate grading shift. Among other things, this gives that the homological thickness of torus link $T_{2k,2kn}$ is (at least) $k(k-1)n+2$, which improves the results obtained in \cite{mojtorus} and \cite{teza}. 

Concerning the first part of the conjecture from \cite{KRS2}, we proved that there exists limit of the rational homology of torus links $T'_{2,2n}$ and $T'_{3,3n}$ when $n$ tends to infinity, and we managed to obtain the limit. Even more, by the method we are using, we calculated explicitly the free part of the homology of all $(3,n)$-torus knots, for every $n\in\N$.\\

\indent In \cite{teza} we have also obtained the analogous stability property of torus knots for Khovanov-Rozansky, or $sl(m)$, homology (see \cite{kovroz}), for every $m\in \N$. However, in the case we are interested in (``cut-off" from above on the homological degrees of torus links), the analogous approach cannot be applied. In addition, we show that for $m$ large enough, there is no such ``cut-off" at all.   \\

The paper iz organized as follows: in Section \ref{not} we give the basic notation and previous results that we need in the paper. In Section \ref{main1} we state the main theorems about the homology of torus links and give their corollaries, while in Section \ref{main2} we compute the center of the ring $H^k$. Section \ref{dokaz} consists of the proofs of the theorems from Section \ref{main1}. In Section \ref{slm} we show that the analogous result for $sl(m)$ for $m>2$ is not valid. Finally, in Section \ref{3n}, we calculate explicitely the free part of the Khovanov homology of $(3,n)$-torus knots, and consequently the stable homology of torus links $T'_{3,3n}$.\\

\textbf{Acknowledgements:} I would like to thank M. Khovanov for suggesting me to study the behaviour of the torus knots at the large homological degrees, for providing me the preprint \cite{KRS2}, and for valuable comments on the early version of the paper. Also, I would like to thank A. Shumakovitch for providing me the results of his computations of the homology of large torus knots.

\section{Notation} \label{not}

For a detailed definition of Khovanov ($sl(2)$) homology of links see e.g. \cite{bn}, \cite{kov}, \cite{mojtorus} or \cite{teza}. Also, we follow the notations from \cite{mojtorus} and \cite{teza} for positive and negative crossing, resolutions, torus knot diagramas, etc.  Here we will just recall a couple of basic facts that  will be used in the paper. 

\indent Let $K$ be a knot (or link) and $D$ its diagram with $n_+$ positive and $n_-$ negative crossings. Then M. Khovanov assigned the graded chain complex $C(D)$, whose homology groups $H(D)$ are bigraded $H^{i,j}(D)$, with $i$ being homological degree and $j$ being second ($q$)-degree. In order to obtain the knot invariant, one has to shift the complex in both degrees:
$$\C(K):=C(D)[-n_-]\{n_+-2n_-\},$$
and the homology groups $\H(K)$ of the complex $\C(K)$ are the knot invariants. More precisely, we have:
\begin{equation}
\label{sift}
\H^{i-n_-,j+n_+-2n_-}(K)=H^{i,j}(D).
\end{equation}

\indent Since for each crossing $c$ of $D$, the complex of the diagram $D$ is the cone of the map between the complexes of the resolutions $D_0$ and $D_1$ of the diagram $D$ at the crossing $c$, we have that for every $i$ and $j$, the following long exact sequence in homology is valid (see e.g. \cite{v2}, \cite{teza}):

{\small{
\begin{equation}\cdots\rightarrow H^{i-1,j-1}(D_1)\to H^{i,j}(D)\to H^{i,j}(D_0)\to H^{i,j-1}(D_1)\to H^{i+1,j}(D)\to \cdots\label{eq1}\end{equation}
}}
 
\indent A knot or a link is a torus knot if it is isotopic to  a knot or a link that can be drawn without any points of intersection on the trivial torus. Every torus link is, up to a mirror image,  determined by two nonnegative integers $p$ and $q$, i.e. it is isotopic to a unique torus knot $T_{p,q}$ which has the diagram $D_{p,q}$ -- the closure of the braid $(\sigma_1\sigma_2\ldots\sigma_{p-1})^q$ -- as a planar projection. In other words, $D_{p,q}$ is the closure of the $p$-strand braid with $q$ full twists. We assume that all strands are oriented upwards, and so that all crossings (in total $(p-1)q$ of them) of $D_{p,q}$ are positive. 

\indent We say that  the crossing $c$ of $D_{p,q}$ is of the type $\sigma_i$, $i<p$, if it corresponds to the generator $\sigma_i$ in the braid word of which $D_{p,q}$ is the closure. For each $i=1,\ldots,p-1$, order the crossings of the type $\sigma_i$ from top to bottom. Then every crossing $c$ of $D_{p,q}$ we can write as the pair $(i,\alpha)$, $i=1,\ldots,p-1$ and $\alpha=1,\ldots,q$, if $c$ is of the type $\sigma_i$ and it is ordered as $\alpha$-th among the crossings of the type $\sigma_i$. \\

\indent Now, consider torus link $D_{2k,2kn}$. Let $c_{2k-1}$ be the crossing $(2k-1,1)$ of the diagram $D_{2k,2kn}$. Denote by $E_{2k,2kn}^1$ and $D_{2k,2kn}^1$ the 1- and 0-resolutions, respectively, of the diagram $D_{2k,2kn}$ at the crossing $c_{2k-1}$. Then from (\ref{eq1}) we obtain
{\small{
\begin{equation*}\cdots\rightarrow H^{i-1,j-1}(E^1_{2k,2kn})\to H^{i,j}(D_{2k,2kn})\to H^{i,j}(D_{2k,2kn}^1)\to H^{i,j-1}(E^1_{2k,2kn})
\to \cdots\end{equation*}
}}
\indent Now, we  continue the process, and resolve the crossing $c_{2k-2}=(2k-2,1)$ of $D_{2k,2kn}^1$ in two possible ways. Denote the diagram obtained by the 1-resolution by $E^2_{2k,2kn}$, and the diagram obtained by the 0-resolution by $D^2_{2k,2kn}$. Then, from the long exact sequence (\ref{eq1}),
we have:
{\small{
\begin{equation*}\cdots\rightarrow H^{i-1,j-1}(E^2_{2k,2kn})\to H^{i,j}(D^1_{2k,2kn})\to H^{i,j}(D_{2k,2kn}^2)\to H^{i,j-1}(E^2_{2k,2kn})\to 
\cdots\end{equation*}
}}
\indent After repeating this process $2k-1$ times (resolving the crossing $c_{2k-l}=(2k-l,1),\,l=1,\ldots,2k-1$, of $D^{l-1}_{2k,2kn}$, obtaining the 1-resolution $E^l_{2k,2kn}$ and 0-resolution $D^l_{2k,2kn}$ and applying the same long exact sequence in homology), we obtain that for every $l=1,\ldots,2k-1$:
{\small{
\begin{equation}\cdots\rightarrow H^{i-1,j-1}(E^l_{2k,2kn})\to H^{i,j}(D^{l-1}_{2k,2kn})\to H^{i,j}(D_{2k,2kn}^l)\to H^{i,j-1}(E^l_{2k,2kn})\to 
\cdots\label{eq3}\end{equation}
}}
Here  $D^0_{2k,2kn}$ denotes $D_{2k,2kn}$, and we obviously have that $D^{2k-1}_{2k,2kn}=D_{2k,2kn-1}$.\\

\indent Also, in the paper we are using some combinatorial results. First of all, we deal with binomial coefficients. If $n$ and $k$ are integers, by $n \choose k$ we define the following expression:
$${n \choose k}=\frac{n!}{k!(n-k)!},$$
if $n$ and $k$ are both nonnegative and $n\ge k$, and zero otherwise. The binomial coefficients obviously satisfy:
$${n\choose k}={n\choose n-k},$$
and 
\begin{equation}
\label{binom}
{n+1 \choose k+1}={n\choose k}+{n\choose k+1}
\end{equation}
which in particular gives
$${2k \choose k}={2k-1 \choose k-1}+{2k-1 \choose k}=2{2k-1\choose k}.$$

Let $k$ be fixed positive integer. By a \textit{zig-zag line} (of length $k$) we mean a graph of a continuous piecewise linear function $f$, defined for $0\le x \le 2k$, such that $f'(x)=\pm 1$ for every $x\in [0,2k]\setminus \Z$ and $f(\Z)\subset \Z$. By a source of a zig-zag line we mean $f(0)$, and by a target we mean $f(2k)$.
Obviuosly, a zig-zag line is uniquely determined by its source and the sequence from $\{-1,1\}^{2k}$ (the sequence of derivatives of a zig-zag line). Hence, the number of zig-zag lines with source $0$ and target $2k-2i$, for $0\le i \le k$ is equal to the number of sequences from $\{-1,1\}^{2k}$, with exactly $2k-i$ entries equal to $+1$ and $i$ entries equal to $-1$, which is equal to $2k \choose i$.
  
Alternatively, we can represent every zig-zag line $f$ by an integral sequence $(a_0,a_1,\ldots,a_{2k})$ with $a_i=f(i)$, $i=0,\ldots,2k$. This gives a bijection between the set of all zig-zag lines and the set of integral sequences of length $2k+1$ such that the difference of the consecutive entries of the sequence is equal to $+1$ or $-1$. So, the number of latter sequences with fixed $a_0$ and $a_{2k}$ is equal to the number of zig-zag lines with source $a_0$ and target $a_{2k}$.


\section{The cut-off in the homology of torus links} \label{main1}

By $T_{2k,2kn}$ we denote the standard $2k$-component torus link with all crossings positive ($(2k-1)2kn$ in total), and by $D=D_{2k,2kn}$ we denote its standard projection (the closure of the braid $(\sigma_1\ldots\sigma_{2k-1})^{2kn})$). By $T'_{2k,2kn}$ we denote the same torus link when exactly $k$ of the components are  with the orientation reversed. Then, the number of positive crossings of $T'_{2k,2kn}$ is equal to $n_+=2k(k-1)n$ and the number of negative crossings is equal to $n_-=2k^2n$. Thus, from (\ref{sift}) we have: 
\begin{eqnarray}
&\H^{i,j+(2k-1)2kn}(T_{2k,2kn})=H^{i,j}(D).\label{x_1}\\
&\H^{i-2k^2n,j-2kn(k+1)}(T'_{2k,2kn})=H^{i,j}(D).\label{x_2}
\end{eqnarray}

Hence, we have:
\begin{equation}
\H^{i,j}(T'_{2k,2kn})=\H^{i+2k^2n,j+6k^2n}(T_{2k,2kn}).
\label{jednacina1}
\end{equation}
In this paper we prove the following result

\begin{theorem}\label{t1}
For every two positive integers $k$ and $n$, we have that $\H^{i}(T_{2k,2kn})$ is trivial for every $i>2k^2n$.
Also, for every $i\in\Z$ the homology group $\H^{i,j}(T_{2k,2kn})$ is trivial for every $j>6k^2n$.
\end{theorem}

The previous theorem, together with (\ref{jednacina1}), gives that the homology of torus link $T'_{2k,2kn}$ is contaned in the lower-left quadrant:

\begin{corollary}
The homology group $\H^{i,j}(T'_{2k,2kn})$ is trivial if $i>0$ or $j>0$.
\end{corollary}

In addition, we obtain the last non-zero homology group of the torus link. Namely, we have:

\begin{theorem}\label{t2}
The rank of the $2k^2n$-th homology group of $T_{2k,2kn}$ is equal to $2k\choose k$, i.e.
$$\rank{\H^{2k^2n}(T_{2k,2kn})}= {2k\choose k},$$ 
and it is torsion-free. More precisely, 
$$\rank{\H^{2k^2n,6k^2n-2i}(T_{2k,2kn})}= {2k\choose k-i} - {2k\choose k-i-1},\quad i=0,\ldots,k,$$
$$\rank{\H^{2k^2n,6k^2n-2i}(T_{2k,2kn})}= 0,\quad \mathrm{for }\quad i<0.$$
\end{theorem}

Again, from (\ref{jednacina1}), we have 

\begin{corollary} \label{cor2}
The rank of the $0$-th homology group of $T'_{2k,2kn}$ is equal to $2k\choose k$, i.e.
$$\rank{\H^{0}(T'_{2k,2kn})}={2k\choose k},$$
 and it is torsion-free. More precisely, 
$$\rank{\H^{0,-2i}(T'_{2k,2kn})}= {2k\choose k-i} - {2k\choose k-i-1},\quad i=0,\ldots,k,$$
$$\rank{\H^{0,-2i}(T'_{2k,2kn})}= 0,\quad \mathrm{for } \quad i<0.$$
\end{corollary}

In other words, the ranks of the groups in the different $q$-gradings in the last non-zero homology group follows the following pattern: for $k=1$ it is $(1,1)$, for $k=2$ it is $(1,3,2)$, for $k=3$ it is $(1,5,9,5)$, etc. \\

\indent Also, from Theorem \ref{t2}, we have that there exists a generator of the homology of $T_{2k,2kn}$, with homological grading equal to $2k^2n$ and $q$-grading equal to $6k^2n-2k$, and hence its $\delta$-grading (see e.g. \cite{mojtorus}) is equal to $6k^2n-2k-4k^2n=2k(nk-1)$. On the other hand, we know that there exists a generator in the homological degree 0, and with the $q$-grading $(2k-1)(2kn-1)+1$ (see e.g. \cite{ras}, \cite{pat}), and hence its $\delta$-grading is equal to $(2k-1)(2kn-1)+1$. Hence, there exist two generators of homology of $T_{2k,2kn}$ whose $\delta$-gradings differ by $(2k-1)(2kn-1)+1-2k(nk-1)=2k(k-1)n+2$, and so we have:
\begin{corollary}
The homological thickness of the torus link $T_{2k,2kn}$ is at least $k(k-1)n+2$.
\end{corollary}

We give the proofs of the Theorems \ref{t1} and \ref{t2} in Section \ref{dokaz}.

\section{Center of the ring $H^k$} \label{main2}

In this section we show that the center of the (graded) ring $H^k$ coincides, up to a grading shift, with the zeroth homology group of torus link $T'_{2k,2kn}$. Namely, we shall prove the following 

\begin{proposition}\label{prop1}
The center of the ring $H^k$ is a free graded abelian group of rank
$${2k \choose i} - {2k \choose i-1}$$
in degree $2i$, for $i=0,\ldots,k$, and zero otherwise.

The total rank of the center is $2k\choose k$.
\end{proposition}

\textbf{Proof:}\\
In the proof we use the representation of the center $Z(H^k)$ obtained in 
\cite{kovspr}. Namely, it was shown that there exists a basis in the center of $H^k$, denoted by $X_I$, which is indexed by all admissible subsets $I\subset \{1,\ldots,2k\}$, where the element $X_I$ is homogeneous of degree $2|I|$.  A subset $I\subset \{1,\ldots,2k\}$ is called \textit{admissible} if  $I \cap \{1,\ldots,m\}$ has at most $m\over 2$ elements for each $m\in\{1,\ldots,2k\}$. For more details see Section 4 of \cite{kovspr}.

\indent Hence, in order to finish the proof, we have to count the number of the admissible subsets with a given cardinality. Obviously, there are no admissible subsets of cardinality strictly bigger than $k$ (just put $m=2k$ in the definition above), and there are no sets of negative cardinality. Thus, there can be admissible subsets only of the cardinality $i$, for $i=0,\ldots,k$. 

\indent Let $i\in\{0,\ldots,k\}$ be fixed. Then there is a bijective correspondence between the set of all subsets of $S=\{1,\ldots,2k\}$ of cardinality $i$, and the set $A_i$ of all integral sequences $(a_0,\ldots,a_{2k})$ of length $2k+1$, such that $a_0=0$, $a_{2k}=2k-2i$, and such that $a_j=a_{j-1}+1$ or $a_j=a_{j-1}-1$, for all $j=1,\ldots,2k$. Indeed, the bijection maps the subset $X\subset S$ to the sequence $a_X=(0,a_1,\ldots,a_{2k})$ given by:
$$a_j=j-2|X\cap\{1,\ldots,j\}|,\quad j=1,\ldots,2k.$$
From the definition it follows that $I$ is an admissible subset if and only if all entries of $a_I$ are nonnegative. We denote the set of such sequences by $A_i^+$, and $B_i=A_i\setminus A_i^+$. Hence, we are left with proving that $|A_i^+|={2k\choose i}-{2k\choose i-1}$.

\indent The cardinality of the set $A_i$ is equal to the number of zig-zag lines of length $2k$ with source $0$ and target $2k-2i$, which is equal to $2k\choose i$ (see end of Section \ref{not}).

\indent Now, let's determine the cardinality of $B_i$. By using the description of the elements of $A_i$ by zig-zag lines, to every element of $a\in B_i$ we can assign bijectively a zig-zag line with the origin $-2$ and the target $2k-2i$. Indeed, since $a\notin A_i^+$, there exists an entry of $a$ equal to $-1$, and by $j$ we denote the smallest index $j$, such that $a_j=-1$. Then, by ``reflecting" the part of the zig-zag line for $x\le j$ with respect to the line $y=-1$, we obtain the zig-zag line  with the source $-2$ and target $2k-2i$, and this correspondence is bijective. Hence the cardinality of the set $B_i$ is equal to the number of zig-zag lines with source $-2$ and target $2k-2i$, which is equal to the number of zig-zag lines with source $0$ and target $2k-2i+2=2k-2(i-1)$, i.e. $2k\choose i-1$. Hence we have:
$$|A_i^+|=|A_i|-|B_i|={2k\choose i}-{2k\choose i-1},$$
which concludes the proof. \kraj
\\

From Proposition \ref{prop1} and Corollary \ref{cor2} we obtain the following result, conjectured in \cite{KRS2}:

\begin{corollary}
For every $k,\,n\in\N$ and $i\in \Z$ the groups $H^{0,-2i}(T'_{2k,2kn})$ and $Z^{2k-2i}(H^k)$ are isomorphic. 

\indent The center $Z(H^k)$ is canonically isomorphic, as a graded abelian group, to $H^0(T'_{2k,2kn})$.
\end{corollary}

\section{Proofs}\label{dokaz}
In this section we prove the theorems announced in Section \ref{main1}. \\

\textbf{Proof of  Theorem \ref{t1}:}\\

\indent From (\ref{x_1}) we have that $\H^i(T_{2k,2kn})$ is trivial if and only if $H^i(D_{2k,2kn})$ is trivial. Thus, we shall prove that $H^i(D_{2k,2kn})$ is trivial for $i>2k^2n$.

In fact, we shall prove even more, that $H^j(D_{2k,2kn}^i)$ is trivial for $j>2k^2n$, and every $i=0,\ldots,2k-2$.\\

We use the induction on $k$. For $k=1$ we have alternating $(2,2n)$ torus links for which the homology is well-known (see e.g. \cite{kov}) and it obviously satisfies the properties of the theorem.\\
\indent Now, suppose that the theorem is valid for $1,\ldots,k-1$, and then we will show that the theorem is valid for $k$. Let $D_{2k,2kn}$ be the standard diagram of the torus link $T_{2k,2kn}$. Then  the diagram $E^1_{2k,2kn}$ has exactly $(4k-2)n-1$ negative crossings, and it is equivalent to the  diagram $D_{2k-2,(2k-2)n} \cup U$, where by $U$ we denoted the unknot. \\
\indent Analogously, the diagram $E^2_{2k,2kn}$ has exactly $(4k-2)n-1$ negative crossings and it is isotopic to the diagram $D_{2k-2,(2k-2)n}$. By proceeding, we obtain that for every 
$i=3,\ldots,2k-1$, 
the  diagram  $E^i_{2k,2kn}$  has exactly $(4k-2)n-1$ negative crossings and that it is isotopic to $D^{i-2}_{2k-2,(2k-2)n}$.\\
\indent Now, by induction hypothesis, we have that for every $i=0,\ldots,2k-3,$  $H^{j}(D^{i}_{2k-2,(2k-2)n})$ is trivial for every $j>2(k-1)^2n$, and so $H^{j}(E^i_{2k,2kn})$ is trivial for every $j>2(k-1)^2n+(4k-2)n-1=2k^2n-1$. Now, by applying long exact sequences (\ref{eq3}), we obtain that  for $i>2k^2n$  the group
$H^{i}(D_{2k,2kn})$ is trivial, if and only if $H^{i}(D^j_{2k,2kn})$ is trivial for every $j=1,\ldots,2k-1$, and in particular, if and only if   
$H^{i}(D_{2k,2kn-1})$ is trivial (and consequently all $H^{i}(D^j_{2k,2kn})$, $j=1,\ldots,2k-2$ are trivial). \\
\indent By repeating the same process, we  obtain that for every $i>k(2kn-1)$, $H^{i}(D_{2k,2kn-1})$ is trivial if and only if   
$H^{i}(D_{2k,2kn-2})$ is trivial. And in general, for every $i>kl$ we have that $H^{i}(D_{2k,l})$ is trivial if and only if $H^i(D_{2k,l-1})$ is trivial.

\indent Finally, since $H^i(D_{2k,2})$ is trivial for $i>2k$, we obtain that for every $i>kl$ we have that $H^{i}(D_{2k,l})$ is trivial, and in particular 
$$H^{i}(D_{2k,2kn})\textrm{ is trivial  for }i>2k^2n,$$
which gives the first part of Theorem 1.

 Analogously we obtain that $H^{i,j}(D_{2k,2kn})$ is trivial  for $j>2k(k+1)n$ and for every $i$. This together with (\ref{x_1}) concludes the proof of Theorem 1.
\kraj\\

\begin{remark}
Completely analogously, we can also obtain that for every $p,q\in \N$, the homology groups $\H^i(T_{p,q})$ are trivial for every $i>pq/2$.
\end{remark}

\textbf{Proof of Theorem \ref{t2}:}\\

\indent Now, we shall concentrate on the last nonzero homology group $$H^{2k^2n}(D_{2k,2kn}).$$ 
In this case we have the following sequence of 
long exact sequences (one for each $i=1,\ldots,2k-1$): 

{\small{
\begin{equation}\cdots\rightarrow H^{2k^2n-1,j-1}(E^i_{2k,2kn})\to H^{2k^2n,j}(D^{i-1}_{2k,2kn})\to H^{2k^2n,j}(D_{2k,2kn}^i)\to  \cdots\label{eqm1}\end{equation}
}}

We shall prove more than stated in Theorem \ref{t2}. Namely, we shall show that 
\begin{equation}\label{gl}
\rank{H^{2k^2n}(D^{i}_{2k,2kn})}=2{2k-1-i\choose k},
\end{equation}
and that $H^{2k^2n}(D^{i}_{2k,2kn})$ is torsion-free for every $i=0,\ldots,2k-1$. \\

\indent  Indeed, we shall prove (\ref{gl}) by induction on $k$. For $k=1$ this is obvious. Now suppose that the statement is true for every $l<k$, and we shall prove that it is valid for $k$. As in the proof of Theorem 1, we have that $H^{2k^2n-1}(E^i_{2k,2kn})=H^{2(k-1)^2n}(D^{i-2}_{2k-2,(2k-2)n})$, for $i=2,\ldots,2k-1$, and 
$H^{2k^2n-1}(E^1_{2k,2kn})=H^{2(k-1)^2n}(D_{2k-2,(2k-2)n}\cup U).$ Hence, 
from the induction hypotheses we have that $H^{2k^2n-1}(E^{i}_{2k,2kn})$ is torsion-free for $i=1,\ldots,2k-1$ and 
$$\rank{H^{2k^2n-1}(E^i_{2k,2kn})}=2{2k-1-i \choose k-1},\textrm{ for } i\ge 2,$$
and   
$$\rank{H^{2k^2n-1}(E^1_{2k,2kn})}= 4{2k-3\choose k-1}=2{2k-2 \choose k-1}.$$

Also, since $H^{2k^2n}(D_{2k,2kn-1})$ is trivial, from (\ref{binom}) and (\ref{eqm1}) we have that 
\begin{eqnarray}
\rank{H^{2k^2n}(D_{2k,2kn})}&\le& \sum_{i=1}^{2k-1} \rank E^i_{2k,2kn} = \label{formula}\\  
&=& 2\sum_{i=k-1}^{2k-2}{i \choose k-1} \\ 
&=&2{2k-1 \choose k}={{2k}\choose k}. 
\end{eqnarray}

On the other hand, the Lee's homology \cite{lee2} of $2k$-component torus link $T_{2k,2kn}$ in homological degree $2k^2n$ is of the rank 
${2k}\choose k$. Since Lee's homology is the $E_{\infty}$-page of the spectral sequence whose $E_2$-page is Khovanov homology (\cite{ras}), we have that 
$$\rank{H^{2k^2n}(D_{2k,2kn})}\ge {{2k}\choose k}.$$
Thus we must have an equality in (\ref{formula}), and so the first map in (\ref{eqm1}) is one-to-one, and the second one is onto. This in turns gives 
$$\rank{H^{2k^2n}(D^j_{2k,2kn})}=2\sum_{i=k-1}^{2k-2-j}{i\choose k-1}=2{2k-1-j \choose k}.$$
In addition, we also obtain the $q$-gradings of the ${2k}\choose k$ 
generators of the last nonzero homology group. Namely, one easily obtain 
$$\rank{\H^{2k^2n,6k^2n-2i}(T_{2k,2kn})}={{2k}\choose k-i}-{2k\choose k-i-1},\quad i=0,\ldots, k$$
and all other homology groups are of the zero rank. Finally, from (\ref{eqm1}), we obviously have that the homology group $H^{2k^2n}(D_{2k,2kn})$ is without torsion, since the homology groups $H^{2k^2n-1,j-1}(E^i_{2k,2kn})$ and $H^{2k^2n,j}(D_{2k,2kn}^i)$ are torsion-free and their ranks sum up to the rank of $H^{2k^2n,j}(D_{2k,2kn}^{i-1})$, which gives Theorem \ref{t2}. Also, note that
\begin{eqnarray*}
\rank{\H^{0,0}(T'_{2k,2kn})}&=&\rank{\H^{2k^2n,6k^2n}(T_{2k,2kn})}={{2k}\choose k}-{2k\choose k-1}=\\
&=&\frac{1}{k+1}{2k\choose k}=C_k,
\end{eqnarray*}
where $C_k$ is the $k$-th Catalan number.\kraj

\section{$sl(m)$ case} \label{slm}

In \cite{teza}, apart from the stability property for Khovanov ($sl(2)$) homology, we have also obtained the analogous results for general Khovanov-Rozansky ($sl(m)$) homology for every $m\in \N$. Namely, if we denote by $H_m(D)$, the $sl(m)$ homology of the diagram $D$, as in \cite{teza} we obtain 
\begin{proposition} 
For every $m,p,q\in\N$ with $p< q$, we have
$$H_m^{i,j}(D_{p,q-1})=H_m^{i,j}(D_{p,q}),\textrm{ for }i<q-1+[q/p](p-2),$$
where by $[x]$ we denoted the largest integer less or equal than $x$.
\end{proposition}

However, in the case that we are interested in, the analog of Theorem \ref{t1} is not valid in general for $m>2$. The basic thing that prevents us from performing the method from \cite{teza} is that the value of the thick edge labelled 3 is nonzero for $m>2$. So, in the $sl(m)$ case we can not reduce the complex from ``above" as we managed from ``below" in  \cite{teza}.

\indent In addition, by using the result from \cite{bojan}, we can easily see that the rank of the (maximal possible) $(2k-1)2kn$-th homology group is at least $m(m-1)\cdots(m-2k+1)$, and hence is nonzero for $m\ge 2k$. For example, in the case of $(4,4n)$-torus links, in $sl(2)$ case, we have proved that the maximal nonzero homological degree is $8n$, for $m=3$ from \cite{bojan} it follows that the rank of the homology group in homological degree $10n$ is at least 6, and for $m\ge 4$ even the homological group in maximal possible degree $12n$ is nonzero. Hence, here the $sl(m)$ homology for $m>2$ significantly differs from the $sl(2)$ one. 

\section{Homology of $(3,n)$-torus knots} \label{3n}

As a simple example how powerful our approach is in the case of torus knots, we calculate the free part of Khovanov homology 
of $(3,n)$-torus knots. Namely we prove the following:

\begin{theorem}\label{t3n}
Poincare polynomial of $(3,3n)$-torus link is given by
\begin{eqnarray*}
P_{3n}(q,t)&=&q^{6n}(q^{-3}+q^{-1}+t^2q+t^3q^5+\\
&&+\left[t^4q^3+t^4q^5+t^5q^7+t^5q^9+t^6q^7+t^7q^{11}\right]\sum_{i=0}^{n-2}{t^{4i}q^{6i}}+\\
&&+t^{4n}q^{6n-3}+3t^{4n}q^{6n-1}+2t^{4n}q^{6n+1}).
\end{eqnarray*}

Poincare polynomial of $(3,3n-1)$-torus knot is given by
\begin{eqnarray*}
P_{3n-1}(q,t)&=& q^{6n-2}(q^{-3}+q^{-1}+t^2q+t^3q^5+\\
&&+\left[ t^4q^3+t^4q^5+t^5q^7+t^5q^9+t^6q^7+t^7q^{11} \right]\sum_{i=0}^{n-2}{t^{4i}q^{6i}}).
\end{eqnarray*}

Poincare polynomial of $(3,3n-2)$-torus knot is given by
\begin{eqnarray*}
P_{3n-2}(q,t)&=&q^{6n-4}(q^{-3}+q^{-1}+t^2q+t^3q^5+\\
&&+\left[t^4q^3+t^4q^5+t^5q^7+t^5 q^9+t^6q^7+t^7q^{11}\right]\sum_{i=0}^{n-2}{t^{4i}q^{6i}}-\\
&& - t^{4n-2}q^{6n-5}-t^{4n-1}q^{6n-1}).
\end{eqnarray*}

\end{theorem}

\indent Like in Section \ref{main1}, we denote by $T'_{2,2n}$ (respectively $T'_{3,3n}$), the torus link which is the same as the standard (positive) torus link $T_{2,2n}$ (respectively $T_{3,3n}$) with one of the components oppositely oriented. Then we have that
\begin{eqnarray*}
&\H^{i,j}(T'_{2,2n})=\H^{i+2n,j+6n}(T_{2,2n}),\\
&\H^{i,j}(T'_{3,3n})=\H^{i+4n,j+12n}(T_{3,3n}).
\end{eqnarray*}

By applying the explicit formula for the homology of torus links $T_{2,2n}$ (already obtained in \cite{kov} and easily obtainable analogously as the result for $(3,n)$-torus knots) and $T_{3,3n}$ (from Theorem \ref{t3n}), we obtain the stable homology of torus links $T'_{2,2n}$ and $T'_{3,3n}$, when $n$ tends to infinity:

\begin{theorem}
There exists limit
$$P_2(q,t)=\lim_{n\to\infty}{\sum_{i,j\in\Z}{t^i q^j \rank \H^{i,j}(T'_{2,2n})}},$$
and it is given by
$$P_2(q,t)=1+q^{-2}+t^{-1}q^{-2}\left(1+t^{-1}q^{-4}\right)\sum_{i=0}^{\infty}t^{-2i}q^{-4i}.$$

Also, there exists limit
$$P_3(q,t)=\lim_{n\to\infty}{\sum_{i,j\in\Z}{t^i q^j \rank \H^{i,j}(T'_{3,3n})}},$$
and it is given by
$$P_3(q,t)=2q+3q^{-1}+q^{-3}+\left(t^{-1}q^{-1}+t^{-3}q^{-3}+t^{-3}q^{-5}\right)(1+t^{-1}q^{-4})\sum_{i=0}^{\infty}t^{-4i}q^{-6i}.$$
\end{theorem}

\textbf{Proof} (of Theorem \ref{t3n}):\\
The main tool is again the long exact sequence (\ref{eq1}). To use it we need the description of the diagrams $E^i_{3,n}$, $i=1,2$, i.e. the number of positive and negative crossing and the knot to which they are isotopic. It is easy to see that we have the following ($U$ denotes the unknot):
\begin{eqnarray}
E^1_{3,3n+3}:&\textrm{ isotopic to }U\cup U,\,n_-=4n+3,\,n_+=2n+2\\
E^2_{3,3n+3}:&\textrm{ isotopic to }U,\,n_-=4n+3,\,n_+=2n+1\\
E^1_{3,3n+2}:&\textrm{ isotopic to }U,\,n_-=4n+2,\,n_+=2n+1\\
E^2_{3,3n+2}:&\textrm{ isotopic to }U,\,n_-=4n+1,\,n_+=2n+1\\
E^1_{3,3n+1}:&\textrm{ isotopic to }U,\,n_-=4n,\,n_+=2n+1 \label{peta}\\
E^2_{3,3n+1}:&\textrm{ isotopic to }U\cup U,\,n_-=4n,\,n_+=2n \label{sesta}
\end{eqnarray} 
Hence we have that, for instance $H^{i,j}(E^2_{3,3n+3})=\Q$ for $i=4n+3$ and $j=6n+5\pm 1$, and zero otherwise.

\indent We prove the formulae from the statement of the theorem by the induction on $n$. For $n=1$ they are true: $(3,1)$-torus knots is isotopic to the unknot, while the values of $(3,2)$ and $(3,3)$-torus knots are already obtained by using the computer programs. Alternatively, one can start from the (trivial) homology of $(3,1)$-torus knots and apply the same process as we do below in the induction step.

\indent  Suppose now that the formulas are true for some $n$. Then the rightmost nonzero homology group of $D_{3,3n}$ is in the homological degree $4n$ and we have:
$$H^{4n,6n-3}(D_{3,3n})=\Q,\quad H^{4n,6n-1}(D_{3,3n})=\Q^3,\quad H^{4n,6n+1}(D_{3,3n})=\Q^2.$$
Then we apply (\ref{eq1}) for the diagram $D^1_{3,3n+1}$ with respect to the crossing $(1,1)$:
{\small{
\begin{equation}0\rightarrow H^{4n,j}(D^1_{3,3n+1})\to H^{4n,j}(D_{3,3n})\to H^{4n,j-1}(E^2_{3,3n+1})\to H^{4n+1,j}(D^1_{3,3n+1})\to 0 \label{eqx1}\end{equation}
}}
Hence, from (\ref{eqx1}) and (\ref{sesta}), we obtain that all homology groups up to the degree $4n-1$ (inclusive) and after the degree $4n+2$ (inclusive) of $D_{3,3n}$ and $D^1_{3,3n+1}$ coincide (the latter ones are all trivial). Also, we have that $H^{4n,6n-3}(D^1_{3,3n+1})=H^{4n+1,6n+3}(D^1_{3,3n+1})=\Q$, and all homology group $H^{4n,j}(D^1_{3,3n+1})$ are trivial for $j<6n-3$ and $j>6n+1$, as well as $H^{4n+1,j}(D^1_{3,3n+1})$ for $j<6n-1$ and $j>6n+3$.\\
\indent Since $D^1_{3,3n+1}$ is two-component link whose linking number (when the two strands are oriented opposite) is equal to $2n$, we have that  there exists a spectral sequence with $E_2$-page being Khovanov homology of $D^1_{3,3n+1}$ and whose $E_{\infty}$-page contains two generators in homological degree 0 and two in homological degree $4n$, with the (Lee's) differential $d_2$ of bidegree $(1,4)$ (see \cite{lee2},\cite{ras}). Hence, the homology group $H^{4n+1,6n-1}(D^1_{3,3n+1})$ must be trivial (since there is no group it can cancel with), and so from (\ref{eqx1}) we have $H^{4n,6n-1}(D^1_{3,3n+1})=\Q^2$. Analogously, we have that the rank of the homology group $H^{4n+1,6n+1}(D^1_{3,3n+1})$ is less then two, and so we have that $\rank H^{4n+1,6n+1}(D^1_{3,3n+1})=$ $\,\rank H^{4n,6n+1}(D^1_{3,3n+1})=x$, with $x$ being equal to 0 or to 1. Hence, up to this last ambiguity, we have determined completely the homology of $D^1_{3,3n+1}$.\\

\indent Now, we proceed by determining the homology of $D_{3,3n+1}$. To this end we apply the long exact sequence (\ref{eq1}) for the diagram $D_{3,3n+1}$, with respect to the crossing $(2,1)$:

{\small{
\begin{eqnarray*}0&\rightarrow& H^{4n,j}(D_{3,3n+1})\to H^{4n,j}(D^1_{3,3n+1})\to H^{4n,j-1}(E^1_{3,3n+1})\to\\
&\to &H^{4n+1,j}(D_{3,3n+1})\to H^{4n+1,j}(D^1_{3,3n+1})\to 0 
\end{eqnarray*}
}}

Thus, from (\ref{peta}) we obtain that all homology groups, except possibly groups at the bidegrees $(4n,6n\pm 1)$, $(4n+1,6n\pm 1)$, of $D_{3,3n+1}$ and $D^1_{3,3n+1}$, coincide. Concerning these four remaining groups, again because of the existence of Lee's differential and since $D_{3,3n+1}$ is a knot (1-component link), we must have that the homology groups $H^{4n+1,6n-1}(D_{3,3n+1})$ and $H^{4n,6n+1}(D_{3,3n+1})$ are trivial, and that 
$$H^{4n,6n-1}(D_{3,3n+1})=\Q\quad \textrm{ and }\quad H^{4n+1,6n+1}(D_{3,3n+1})=\Q.$$ Thus, we have obtained the required formula for $(3,3n+1)$-torus knots.

\indent The process for obtaining the two remaining homologies (of $D_{3,3n+2}$ and of $D_{3,3n+3}$) is completely analogous to the previously described. 
One easily obtains that the homology of $D_{3,3n+2}$ coincides with the homology of $D_{3,3n+1}$ except at two bidegrees where we have:
$$H^{4n+2,6n+1}(D_{3,3n+2})=H^{4n+3,6n+5}(D_{3,3n+2})=\Q.$$
Finally, for $D_{3,3n+3}$ one obtains that its homology coincides with the homology of $D_{3,3n+2}$ except at three bidegrees where it is as follows:
$$H^{4n+4,6n+3}(D_{3,3n+3})=\Q$$
$$H^{4n+4,6n+5}(D_{3,3n+3})=\Q^3$$ $$H^{4n+4,6n+7}(D_{3,3n+3})=\Q^2.$$
This concludes the proof. \kraj\\

\begin{remark} The result of Theorem \ref{t3n} was also independently obtained by P. Turner, \cite{tur}. 
\end{remark}

\end{document}